\documentstyle{amsppt}
\magnification=\magstep1
\topmatter
\title
A conjecture for the sixth power moment of the Riemann zeta-function
\endtitle
\author
J. B. Conrey 
\\
A. Ghosh 
\endauthor
\thanks Research of the first author supported by the American Institute
of Mathematics. Research of both authors
supported in part by
 a grant from the NSF.
\endthanks
\address American Institute of Mathematics, 360 Portage Avenue, Palo Alto,
California 94306
\endaddress
\address Department of Mathematics, Oklahoma State University,
Stillwater, Oklahoma 74078 
\endaddress
\endtopmatter
\define \z{\zeta}

\bigskip

In 1918 Hardy and Littlewood [2] proved that
$$\int_1^T|\z(1/2+it)|^2~dt \sim T\log T $$
and in 1926 Ingham [4] showed that 
$$ \int_1^T |\z(1/2+it)|^4~dt \sim \frac{1}{2 \pi ^2} T \log^4 T .$$
In general, it is conjectured that if $k>0$, then there exists a $c_k>0$
such that 
$$ \int_1^T |\z(1/2+it)|^{2k}~dt \sim c_k T \log^{k^2} T .$$
No value has been suggested for $c_k$ if $k$ is different from 0,1, or 2.

In this paper, we present evidence to support the

\proclaim{Conjecture} As $T\to \infty$,
$$  \int_1^T |\z(1/2+it)|^6~dt \sim \frac{42}{9!} 
\prod_p\left\{\left(1-\frac1p\right)^{4}\left(1+\frac4p+\frac1{p^2}\right)\right\}T\log^9T
 .$$
\endproclaim

\medskip

\centerline{SKETCH OF BASIC ARGUMENT}

\medskip

We recall the functional equation of the zeta-function:
$$\zeta(s)=\chi(s)\zeta(1-s)$$
where $$\chi(1-s)=\chi(s)^{-1}=2(2\pi)^{-s}\Gamma(s)\cos \frac{\pi s}2.$$
We also require an ``approximate'' functional equation for $\zeta(s)^2$:
$$\zeta(s)^2=D(s)+\chi(s)^2D(1-s)$$
where 
$$D(s)=\sum_{n\le \frac{|t|}{2\pi}}\frac{d(n_)}{n^s}+E(s)$$
with $d(n)$ the usual divisor function and with $E(s)$ a suitable
error-term. The estimate
$$E(1/2+it)\ll \log(2+|t|)$$ 
is known to hold (see [3]). However, for the purposes of this paper we
will not be concerned with $E(s)$.

The beginnings of our argument are as follows.

$$\split
\int_0^T|\zeta(1/2+it)|^6~dt&=\frac1i\int_\frac 12^{\frac12+iT}
\zeta(s)^3\zeta(1-s)^3~ds\\
&=
\frac1i\int_\frac 12^{\frac12+iT}
\chi(1-s)\zeta(s)^4\zeta(1-s)^2~ds\\
&=
\frac1i\int_\frac 12^{\frac12+iT}
\chi(1-s)\zeta(s)^4\left(D(1-s)+\chi(1-s)^2D(s)\right)~ds\\
&=
\frac1i\int_\frac 12^{\frac12+iT}
\chi(1-s)\zeta(s)^4D(1-s)~ds\\
&\qquad +\frac1i\int_\frac 12^{\frac12+iT}
\chi(1-s)^3\zeta(s)^4D(s)~ds\\
&=I_1+I_2,
\endsplit
$$
say.
Now, 
$$
\split
I_2&=\frac1i\int_\frac 12^{\frac12+iT}
\chi(1-s)^3\chi(s)^4\zeta(1-s)^4D(s)~ds\\
&=
\frac1i\int_\frac 12^{\frac12+iT}
\chi(s)\zeta(1-s)^4D(s)~ds\\
&=
\frac1i\int_{\frac 12-iT}^{\frac12}
\chi(1-s)\zeta(s)^4D(1-s)~ds\\
&=\overline{I_1}.
\endsplit
$$

Thus, 

$$  \int_1^T |\z(1/2+it)|^6~dt=2\Re \frac1i\int_\frac 12^{\frac12+iT}
\chi(1-s)\zeta(s)^4D(1-s)~ds.$$

\medskip

\centerline{THEOREM FROM ``MEAN-VALUES III''}

\medskip

We appeal to Theorem 2 of Conrey - Ghosh [1] to evalute
this integral.
We first set up some notation
so that we can state a special case of that Theorem.
Define 
 $D_N(s,P)$ by
$$
D_N(s,P)=\sum_{n\le N}\frac{d(n)}{n^s}P(\log n/\log N)
$$
where $P$ is any real polynomial.
Let 
$$\split
K_{N}(T)&=\int_1^T|\zeta(1/2+it)|^2
\zeta(1/2+it)^2{D_N(1/2-it,P)}~dt
\\
&=\frac1i\int_\frac 12^{\frac12+iT}
\chi(1-s)\zeta(s)^4D_N(1-s,P)~ds.
\endsplit
$$

\proclaim{Theorem}If 
$N=T^\theta$ with $0<\theta<1/2$, then
$$ K_{N}(T) \sim
 T(\log N)^9 \frac{a_3}
{720\theta^3}
\int_0^1P(\alpha)\alpha^{5}{}_2F_1(-2,-3,6,-\alpha\theta)~d\alpha
$$
as $T \to \infty$ where ${}_2F_1$ is the usual hypergeometric function
and
$$a_3=
\prod_p\left\{\left(1-\frac1p\right)^{4}\left(1+\frac4p+\frac1{p^2}\right)\right\}
.$$
\endproclaim
The hypergeometric function simplifies to
$$1-\alpha\theta+(\alpha\theta)^2/7.$$
We are interested in the case where $P=1$.
We find in this case that 
$$
\split
K_N(T)&\sim
\frac{a_3}{\theta^33!5!}\int_0^1\alpha^5\left(1-\alpha\theta+\frac{(\alpha\theta)^2}{7}\right)~d\alpha
T\log^9N
\\
&=\frac{\theta^6a_3}{720}\left(\frac16-\frac{\theta}{7}+\frac{\theta^2}{56}\right)T\log^9T.
\endsplit
$$

This formula is valid for $\theta<1/2$. To apply it 
to our formula from the last section, we would need it to hold for 
$\theta=1$. In our paper [1], we expressed the belief that
the Theorem from which the above is taken is actually valid for all 
$\theta \le 1$.

If we assume that we can take $\theta=1$ in this Theorem,
then we are immediately led to
$$\int_0^T|\zeta(1/2+it)|^6~dt \sim 42 \frac {a_3}{9!} T\log^9T .$$

\medskip

\centerline{SKETCH OF ANOTHER METHOD}

\medskip

We remark that we can arrive at the same conclusion by another method,
which we briefly sketch.
 
With $s=1/2+it$, we have
$$\split
\int_0^T  |\zeta(1/2+it)|^6~dt&=\int_0^T|\zeta(s)|^2\left|
\zeta(s)^2\right|^2~dt
\\
&=\int_0^T|\zeta(s)|^2\left|D(s)+\chi(s)^2D(1-s)\right|^2~dt\\
&=2\int_0^T|\zeta(s)|^2|D(s)|^2~dt+2\Re
\int_0^T|\zeta(s)|^2\chi(1-s)^2D(s)^2~dt
\endsplit
$$
since $|\chi(1/2-it)|=1$.

To evaluate the first integral here we appeal to a
special case of Theorem 1 of Conrey - Ghosh [1].
\proclaim{Theorem} Let
$$
J_{N}(T)=\int_1^T|\zeta(1/2+it)|^2|D_N(1/2+it,P)|^2dt
$$
If $N=T^\theta$ for some $\theta$ with $0<\theta<1/2$,
and if $P$ is a real polynomial, then 
$$
J_{N}(T)
 \sim T( \log N)^{9} \frac{a_3}{24}
\int_0^1\alpha^{3}\left(\frac1\theta h'(\alpha)^2+4 h(\alpha)h'(\alpha)
\right)~d\alpha
$$ as $T \to \infty$ 
where
$$
h(\alpha)=\int_{\alpha}^1(\beta-\alpha)^2 P(\beta)~d\beta.
$$
\endproclaim

Again this Theorem can be proven for $\theta<1/2$ and again
we expressed the belief in [1] that it actually holds true
for $\theta\le 1.$ Assuming the formula for $\theta=1$
leads to 
$$2\int_0^T|\zeta(s)|^2|D(s)|^2~dt\sim 28 \frac {a_3}{9!} T\log^9T ~dt.$$

To handle the second integral we appeal again to the approximate
functional equation for $\zeta(s)^2$.
We find that the second term above is 
$$
\split
&=2\Re \int_0^T\chi(1-s)^3\zeta(s)^2D(s)^2~ds
\\
&=2\Re \int_0^T\chi(1-s)^3\left(D(s)+\chi(s)^2D(1-s)\right)D(s)^2~ds
\\
&=2\Re \int_0^T\chi(1-s)^3D(s)^3~ds+
2\Re \int_0^T\chi(1-s)D(1-s)D(s)^2~ds.
\endsplit
$$
The first integral here is not expected to contribute to the main term,
essentially because
$$\chi(1/2-it)^3=\exp\left(3it\log\frac{t}{2\pi e}\right)$$
is ``spinning'' too fast. 
To evaluate the second integral, we proceed as in the proof
in [1] of Theorem 2.  If the $D(1-s)$ were replaced by
$D_N(1-s,1)$ with $N=T^\theta$ and  $\theta<1/2$ then, 
in a way similar to the proof of Theorem 2, we could obtain
an asymptotic evaluation. Again we assume that this asymptotic
evaluation is actually correct for all $\theta\le 1$.
In this way we obtain
$$2\Re \int_0^T\chi(1-s)D(1-s)D(s)^2~ds
\sim 14\frac{a_3}{9!}T\log^9T.$$

This argument again leads to the conjecture of this paper.

\medskip

\centerline{FINAL REMARKS}

\medskip

Conrey and Gonek, in work in progress, have arrived 
at exactly the same conjecture using a third method. Their method
is to consider the asymptotic behavior of ``long'' Dirichlet polynomials,
based on techniques developed by Goldston and Gonek. Their
method supposes that 
asymptotic formulae exist for sums
$$\sum_{n\le x} d_3(n)d_3(n+h),$$
where $d_3$ can be defined by
$$\zeta(s)^3=\sum_{n=1}^\infty\frac{d_3(n)}{n^s},$$
and that the asymptotic formulae for these sums have 
smooth main-terms and error terms which are bounded 
on average over $h$ by $x^{1/2+\epsilon}$ for 
$h$ up to $x^{1/2+\epsilon}$.

Finally, we mention that another possibility for testing our conjecture
would be through the method
recently developed in the thesis of Jose Gaggero Jara
(under the direction of S. M. Gonek at the University of Rochester).
In that work Jara develops an asymptotic formula for
$$\int_0^T|\zeta(1/2+it)|^4\left|\sum_{n=1}^N a_n n^{it}\right|^2~dt$$
for arbitrary positive coefficients $a_n$ provided that
$N=T^\theta$ with $\theta <4/589$.  It is probably the case that
the formula should actually hold for all $\theta<1/2$ and 
possibly even for all $\theta <1$. With $\theta=1/2$ and
$a_n=1$, the result should give one-half of the sixth moment.
With $\theta=1$ and $a_n=1$ it should give all of the sixth moment.

Generally, by taking an approximate functional equation for
$\zeta(1/2+it)$ with ``uneven'' lengths $t^\theta$ and $t^{1-\theta}$
the sum of the results of Jara's theorem with $a_n=1$ and
 $N=T^\theta$ and $N=T^{1-\theta}$
(for any $\theta < 1$) 
should also give the sixth moment.

\enddocument